\documentstyle[11pt]{article}
\normalbaselines
\begin{document}
\bibliographystyle{unsrt}

\def\bea*{\begin{eqnarray*}}
\def\eea*{\end{eqnarray*}}
\def\ba{\begin{array}}
\def\ea{\end{array}}
\count1=1
\def\be{\ifnum \count1=0 $$ \else \begin{equation}\fi}
\def\ee{\ifnum\count1=0 $$ \else \end{equation}\fi}
\def\ele(#1){\ifnum\count1=0 \eqno({\bf #1}) $$ \else \label{#1}\end{equation}\fi}
\def\req(#1){\ifnum\count1=0 {\bf #1}\else \ref{#1}\fi}
\def\bea(#1){\ifnum \count1=0   $$ \begin{array}{#1}
\else \begin{equation} \begin{array}{#1} \fi}
\def\eea{\ifnum \count1=0 \end{array} $$
\else  \end{array}\end{equation}\fi}
\def\elea(#1){\ifnum \count1=0 \end{array}\label{#1}\eqno({\bf #1}) $$
\else\end{array}\label{#1}\end{equation}\fi}
\def\cit(#1){
\ifnum\count1=0 {\bf #1} \cite{#1} \else 
\cite{#1}\fi}
\def\bibit(#1){\ifnum\count1=0 \bibitem{#1} [#1    ] \else \bibitem{#1}\fi}
\def\ds{\displaystyle}
\def\hb{\hfill\break}
\def\comment#1{\hb {***** {\em #1} *****}\hb }

\newcommand{\TZ}{\hbox{\bf T}}
\newcommand{\MZ}{\hbox{\bf M}}
\newcommand{\ZZ}{\hbox{\bf Z}}
\newcommand{\NZ}{\hbox{\bf N}}
\newcommand{\RZ}{\hbox{\bf R}}
\newcommand{\CZ}{\,\hbox{\bf C}}
\newcommand{\PZ}{\hbox{\bf P}}
\newcommand{\QZ}{\hbox{\bf Q}}
\newcommand{\HZ}{\hbox{\bf H}}
\newcommand{\EZ}{\hbox{\bf E}}
\newcommand{\GZ}{\,\hbox{\bf G}}
\newcommand{\DZ}{\, \hbox{\bf D}}
\newcommand{\FZ}{\hbox{\bf F}}
\newcommand{\KZ}{\hbox{\bf K}}

\newtheorem{theorem}{Theorem}
\newtheorem{lemma}{Lemma}
\newtheorem{proposition}{Proposition}
\newtheorem{corollary}{Corollary}

\vbox{\vspace{38mm}}
\begin{center}
{\LARGE \bf Intermediate Jacobian and Some Arithmetic \\[2mm]
Properties of Kummer-surface-type CY 3-folds }\\[10 mm]  
Shi-shyr Roan
\\{\it Institute of Mathematics, Academia Sinica \\ 
Taipei , Taiwan \\ (e-mail: maroan@gate.sinica.edu.tw)} 
\\[35mm]
\end{center}

\begin{abstract}
In this article, we examine the arithmetic aspect of the Kummer-surface-type CY 3-folds $\widehat{T/G}$, characterized by the crepant resolution of 3-torus-orbifold $T/G$ with only isolated singularities. Up to isomorphisms, there are only two such space $\widehat{T/G}$ with $|G|=3, 7$, and both $T$ carrying the structure of triple-product structure of a CM elliptic curve. The (Griffiths) intermediate Jacobians of these $\widehat{T/G}$ are identified explicitly as the corresponding elliptic curve appeared in the structure of $T$. We further provide the $\QZ$-structure of $\widehat{T/G}$ and verify their modularity property. 
\end{abstract}

\par \vspace{5mm} \noindent

\vfill
\eject

\section{Introduction}
In algebraic geometry study of projective manifolds with trivial canonical bundle, orbifold geometry of (complex) $n$-torus $T$ quotiented by a (finite) special automorphism group $G$ naturally appears in the context due to the characteristic role of Kummer surfaces in the class of K3 surfaces for $n=2$, as one tends to follow the path of Kummer-surface to construct  higher dimensional Calabi-Yau (CY) varieties. 
In doing so, not only arises the problem of existence of crepant resolution of $T/G$, a topic which deserves the study by its own right \cite{Rtop}, but also emerges the structure of tori $T$  with certain constraints (for a special class of group $G$, see \cite{R89}). For $n=3$, two rigid 3-folds $\widehat{T/G}$ were constructed along this line in the early study of CY spaces \cite{RY}, with $(T, G)$ isomorphic to $(E_\omega^3, \ZZ_3)$ or $(A(\mu), \ZZ_7)$ where $E_\omega^3$ is the triple-product of the elliptic curve $E_\omega$ with an order 3 automorphism $m_\omega$, and $A(\mu)$ is an abelian variety of CM type associated to the cyclotomic field of degree 7 (for the precise description, see (\req(Z7)) of this paper). The geometric properties of $T/G$ revealed in \cite{RY}, which are associated to the $G$-fixed points of $T$, strongly suggested the uniqueness of those two Kummer-surface-type CY 3-folds. This lead to the full investigation of all 3-tori $T$ carrying a special automorphism group $G$ which contains an order 7 element in \cite{R7}, where by number-theoretic methods  on abelian varieties with complex multiplications developed in \cite{ST} and known facts on  class numbers of cyclotomic fields with a small degree, the classification of all such $(T, G)$ was obtained. In particular, these two 3-folds $\widehat{T/G}$ constructed in \cite{RY} are indeed the only Kummer-surface-type CY 3-folds, i.e., CY 3-folds $\widehat{T/G}$ with ${\rm Sing}(T/G) < \infty$. On the other hand, in the classification of minimal threefolds with zero Kodaira-dimension, there stood out a class of CY threefolds of type ${\rm III}_0$ \cite{O93}. Based on a structure theorem in \cite{SW}, these two Kummer-surface-type CY 3-folds $\widehat{T/G}$ reappeared in the context of birational geometry, with characterization as the simply-connected ${\rm III}_0$-type CY threefolds \cite{O}. A pedagogical account of these two special CY 3-folds $\widehat{T/G}$ can be found in \cite{R03} \S 2, where some essential number-theoretic facts in \cite{R7} relevant to the justification of uniqueness of two Kummer-surface-type CY 3-folds $\widehat{T/G}$ were put in the appendix for the reason not only due to the self-contained nature, but further for indication of some interesting questions in connection with 3-tori carrying certain non-special automorphism, though not in the context of CY spaces, that remain unanswered, such as the clearer relationship between a class of 3-tori and the spectral curves appeared in a renowned 2-dimensional statistical solvable model: the chiral Potts $N$-state model  \cite{B91} for $N=4$. Nevertheless, in the context of Kummer-surface-type CY 3-folds, their rigid property should be considered as a defect from the geometric aspect in regard to the classification of projective 3-folds. However, the specialty of the spaces $\widehat{T/G}$ and their simple geometrical construction through number fields  appear to be of great interest in searching their arithmetic properties from the number-theoretic point of view. That is the general conviction of our motivation to explore the possible arithmetic structures of these two specific CY 3-folds, especially those relevant to the now well-publicized "modularity property" of rigid CY 3-folds \cite{Yui}. The objective of this paper is to identify the explicit structure of their (Griffiths) intermediate Jacobian, provide a $\QZ$-model and verify the modularity of these two spaces $\widehat{T/G}$.

In this paper, we fix the third and seventh root of unity
$$
\omega := {\rm e}^{2\pi {\rm i}/3} \ , \ \ \ \mu := {\rm e}^{2\pi {\rm i}/7} \ .
$$
and set $
\eta := \mu + \mu^2  + \mu^4 = \frac{-1 + \sqrt{7} {\rm i}}{2} $,
which satisfies the relation
\be
\eta^2 + \eta + 2 = 0 \ .
\ele(eta)
The Kummer-surface-type CY 3-folds $\widehat{T/G}$ are described as follows: for $|G| =3$, 
$(T, \ZZ_3) \simeq (E_\omega^3 , \langle m_\omega^3 \rangle )$, and for $|G| =7$, there are three models to present $(T, G)$ \cite{R7, R03},
\be
(T, \ZZ_7) ~ \simeq ({\rm Jac}(K) , \ \kappa_7 ) ~ \stackrel{(i)}{\simeq} ( A( \mu ) , \langle m_\mu \rangle )  ~ \stackrel{(ii)}{\simeq} (E_\eta^3  , \ \langle g_\eta \rangle ) \ , 
\ele(Z7)
where $\kappa_7$ is the special automorphism group of the Jacobian of the Klein quartic 
$K$ \cite{Kl},
\be
K : \ \ Z_1^3 Z_2 + Z_2^3 Z_3 + Z_3^3 Z_1 = 0 \ , \ \ [ Z_1 , Z_2 , Z_3 ] \in \PZ^2 \ ,
\ele(Kqrt)
induced by an order 7 automorphism of $K$, $A( \mu )$ is the abelian variety associated to the CM field $(\QZ (\mu), \{ \iota_4, \iota_2, \iota_1 \})$ for $\iota_k( \mu)= \mu^k$ with $m_{\mu}$ the special automorphism of $A( \mu )$ induced by the $\mu$-multiplication\footnote{The $A( \mu )$ here was denoted by $A( \QZ(\mu) )$ in \cite{R03}. However for convenience of discussions in this paper, the uniformizing coordinates of the torus differ than that in \cite{R03} by a permutation of indices: $(Z_1, Z_2, Z_3)$ of $A( \mu )$ in this paper equals to $(z_3, z_2, z_1)$ for the coordinates $(z_1, z_2, z_3)$ of $A( \QZ(\mu) )$ in \cite{R03}.}, and $g_{\eta}$ is the following special automorphism of the triple-product $E_\eta^3$ of 1-torus $E_\eta \ (:= \CZ/(\ZZ + \ZZ \eta ) )$,
\be
g_\eta := \left( \begin{array}{ccc} 0 & 0 & 1 \\
1 & 0 & \eta + 1 \\
0 & 1 & \eta \end{array}
\right).
\ele(eta7)
(Note that the eigenvalues of $g_\eta$ are $\mu^4, \mu^2, \mu$.) Even though one can verify the relations in (\req(Z7)) by the uniqueness property of $(T, \ZZ_7)$, it will be clear that an explicit description of isomorphisms between the different models in (\req(Z7)) is important for the study of arithmetic properties of the CY 3-fold, as indicated in a sequel of the present article. For both two CY 3-folds $\widehat{T/G}$ in our consideration, the 3-torus $T$ carries the structure of triple-product of an elliptic curve $E$, from which the intermediate Jacobian of $\widehat{T/G}$ will be derived, and stated in Theorem \ref{thm:InJ} of this paper. The discussion of $\QZ$-structure of $\widehat{T/G}$ will be based on the $\QZ$-equations of $E_\omega$ and $K$ compatible with the group action of $G$. The results are given in Theorem \ref{thm:3eqn} and \ref{thm:7eqn}. Using the triple-product structure $E^3$ for $T$, one derives the modularity of $\widehat{T/G}$ from the modular property of $E$ and cohomology description of $\widehat{T/G}$.

The remainder of this article is organized as follows.  In Sec. 2, we first give a short review of some basic facts on Klein quartic relevant to the present work by following the beautiful expository paper of N. Elkies on the subject \cite{El}, then establish the direct connections between the three models in (\req(Z7)) for later use. In Sec. 3, we derive the structure of the intermediate Jacobian of $\widehat{T/G}$. By the structure of  triple-product  of 1-torus $E$ for $T$, the intermediate Jacobian of $\widehat{T/G}$ is given by the corresponding 1-torus $E$; in both cases, the intermediate Jacobian is a CM elliptic curve defined over $\QZ$. In Sec. 4, we study the rational structure of $T/G$ and $\widehat{T/G}$, and derive their $\QZ$-models. In Sec. 5, we verify the modularity of the rigid CY 3-folds $\widehat{T/G}$.
We end with the concluding remarks in Sec. 6. 

{\bf Convention}. 
To present our work, we prepare some notations. In this paper, 
$\ZZ, \RZ, \CZ$ will denote the ring of integers, real, complex numbers
respectively, $\HZ_+ = \{ \tau \in \CZ \ | \ \Im \tau > 0 \}$ the upper-half plane, and $\ZZ_N=
\ZZ/N\ZZ$, ${\rm i} = \sqrt{-1}$. For a positive integer $n$, $\CZ^n$ is the $n$-dimensional vector space over $\CZ$ consisting of column vectors. For a (complex) $n$-torus $T \ (= \CZ^n/L)$, a group of  Lie-automorphisms of $T$ persevering its holomorphic $n$-form will be called a special automorphism group of $T$.
For $\tau \in \HZ_+$, $E_{\tau}$ will always denote the elliptic curve $\CZ/ (\ZZ + \ZZ \tau)$ with the modular $\tau$, and $E_\tau^n$ the $n$-product of $E_\tau$. When $\tau$ satisfies a monic quadratic integral polynomial equation, or equivalently, $(X- \tau) (X - \overline{\tau}) \in \ZZ[X]$,  
$\ZZ [\tau]$ can be regarded as an endomorphism ring of $E_\tau$, induced by  scalar-multiplications of $\ZZ [\tau]$ on $\CZ$. Accordingly, one may consider ${\rm GL}_n(\ZZ [\tau])$ as an automorphism group of the $n$-torus $E_\tau^n$.


\section{Jacobian of Klein Quartic}
For later use in this paper, we are going to derive the direct links among various models in (\req(Z7)). Moreover, we will show that the symmetries attached to isomorphisms in (\req(Z7)) can further be extended to the non-cyclic group of order 21.

It is known that the Klein quartic $K$ is the genus 3 curve with the automorphism group  isomorphic to the (unique) order 168 simple group $H_{168}$ ( $\simeq {\rm PSL}_2(\FZ_7) \simeq {\rm SL}_3(\FZ_2) $). We consider the Klein model of $H_{168}$ in ${\rm SL}_3(\CZ)$ generated by the following matrices:
\be
g = \left( \begin{array}{lll} \mu^4 & 0 & 0 \\
0 & \mu^2 & 0 \\
0 & 0 & \mu \end{array} \right) \ , \ \ \ 
h = \left( \begin{array}{lll} 0 & 1 & 0 \\
0 & 0 & 1 \\
1 & 0 & 0 \end{array} \right) \ , \ \ \ r = \frac{-1}{\sqrt{7} {\rm i}}\left( \begin{array}{lll} \mu-\mu^6 & \mu^2-\mu^5 & \mu^4-\mu^3 \\
\mu^2-\mu^5 & \mu^4-\mu^3 & \mu-\mu^6 \\
\mu^4-\mu^3 & \mu-\mu^6 & \mu^2-\mu^5 \end{array} \right) \ . 
\ele(ghr)
The above matrices satisfy the following relations:
\be
g^7 = h^3 = r^2 = (r g^4)^4 = (rh)^2 = 1 , \ h^{-1}gh = g^2  . 
\ele(21168)
By a well-known fact in finite group theory \cite{Bl, Bu}, the above relations in $g, h$ characterize $\langle g, h \rangle$ as the unique non-cyclic group of order 21 (up to isomorphisms), and all the relations in (\req(21168)) give the characterization of $H_{168}$. By the matrices in (\req(ghr)), $H_{168}$ acts on $\CZ^3$ with the coordinates $(Z_1, Z_2, Z_3)^t$. The invariant ring $\CZ[Z_1, Z_2, Z_3]^{H_{168}}$ is generated by three algebraic independent homogenous polynomials $\Phi_4, \Phi_6, \Phi_{14}$ of degrees 4, 6, 14 and a fourth one $\Phi_{21}$ whose square is a polynomial of the first three \cite{Kl} (or see \cite{El} \S 1). The polynomial $\Phi_4$ is the one in the equation (\req(Kqrt)) of Klein quartic $K$. The CM field $(\QZ (\mu), \{ \iota_4, \iota_2, \iota_1 \})$, which is compatible with the matrix expression of $g$, gives rise to the $\RZ$-isomorphism $u : \QZ(\mu) \otimes \RZ \longrightarrow \CZ^3$, where $u (x) = (\iota_4(x) , \iota_2 (x), \iota_1 (x) )^t$ for $x \in \QZ(\mu)$. The abelian variety $A( \mu )$ is the 3-torus defined by
$$
A( \mu ) = \CZ^3 / u( \ZZ[\mu] ) \, \, \ \ \ \   \big( \simeq (\QZ(\mu) \otimes \RZ)/\ZZ[\mu]  \, \ {\rm as \ real \ tori} \big) \ . 
$$
Indeed, $\ZZ[\mu]$ is the ring of integers of $\QZ[\mu]$, and the lattice $u(\ZZ[\mu])$ in $\CZ^3$ has the following integral basis:
$$
\bigg( u(1), u(\mu), \cdots , u(\mu^5) \bigg)= \left( \begin{array}{cccccc}
1 & \mu^4 & \mu & \mu^5 & \mu^2 & \mu^6 \\
1 & \mu^2 & \mu^4 & \mu^6 & \mu & \mu^3 \\
1 & \mu & \mu^2 & \mu^3 & \mu^4 & \mu^5 \\
\end{array} \right) \ .
$$
The special automorphism $m_{\mu}$ of  $A( \mu )$ corresponds to the one induced by the matrix $g$ on the universal cover $\CZ^3$ of $A( \mu )$. Furthermore, the matrix $h$ induces the a $\ZZ$-isomorphism of the lattice $u( \ZZ[\mu] )$. Hence $\langle g, h \rangle$ gives rise to the order 21 special automorphism group of $A( \mu )$. But the discussion can not be carried further to $H_{168}$ as the lattice $u( \ZZ[\mu] )$ is not invariant under the matrix $r$. 
Note that the zeros, $\Phi_4 = 0$, in $\CZ^3$ do not produce a curve in $A( \mu )$ by $\Phi_4 (u(1))=4$. The $m_\mu$-fixed point set $A( \mu )^{m_\mu }$ is the order 7 (additive) subgroup of $A( \mu )$ with the generator $\alpha: = \frac{1}{7} \sum_{j=0}^5 u(\mu^j)$, and is stable under the special automorphism $h$ of $A( \mu )$ which sends $\alpha$ to $2 \alpha$. Hence $A( \mu )^{m_\mu }$ consists of  three $\langle h \rangle$-orbits, $\{ 0 \}$, $\{ \alpha , 2 \alpha , 4 \alpha \}$, $\{ 3\alpha , 6 \alpha , 5 \alpha \}$.

We now consider $H_{168}$ as a subgroup of ${\rm PGL}_2(\CZ)$ through the representation (\req(ghr)), acting on the Klein quartic $K$ as the automorphism group. The quotient of $K$ by the subgroup $\langle g \rangle $ is a rational curve, by which $K$ is a cyclic cover of $\PZ^1$ of degree 7 birational to the plane curve,  
\be
Y^7 = X^2 (X-1) \ , \ \ (X, Y ) \in \CZ^2 \ ,
\ele(Kaff)
where $g$ corresponds to the map, $(X, Y) \mapsto ( X, \mu Y)$ (see, e.g. \cite{R03}). Note that the $g$-fixed point set $K^g$ in $K$ consists of three elements in (\req(Kaff)) with the $X$-values $0, 1, \infty$, denoted by $q_0, q_1, q_\infty$ respectively, which as elements in (\req(Kqrt)), have the corresponding  homogenous coordinates given by $[Z_1, Z_2, Z_3]=[1, 0, 0], [0,0, 1] , [0, 1, 0]$. Let $t_j$ be the local smooth coordinate of $K$ centered at $q_j \in K^g$. The automorphism $g$ has the local description 
\be
g: t_j \mapsto u_j t_j, \ \ \ \ \ u_j=u_j(t_j) : {\rm a \ unit \ near} \ q_j \ \ {\rm  with} \ u_0(0)= \mu^4, \ u_1(0)= \mu, \ u_\infty(0)= \mu^2. 
\ele(gloc) 
Denote $\varphi_1= -dX/ Y^3$ , $ \varphi_2= XdX/ Y^5 $, $\varphi_3= XdX/ Y^6$. Then $\varphi_j$'s form a basis of abelian differentials (of the first kind)\footnote{Here we choose $-dX/ Y^3$ instead of $dx/ y^3$ in \cite{R03} as a basis element for the purpose of later discussion on the $h$-action of $K$.} for $K$, by which the $g$-induced automorphism of Jacobian of $K$, 
$$
 {\rm Jac}(K) = \Gamma (K, \Omega)^*/H_1(K, \ZZ) \ ,
$$
is compatible the matrix form of $g$ in (\req(ghr)). For simplicity, the $g$-induced automorphism of ${\rm Jac}(K)$ will be denoted again by $g$ if no confusion could arise.  Then $g$ is a generator of the group $\kappa_7$ in (\req(Z7)). The explicit relations between 
$X, Y$ and $Z_j$'s are given by \cite{El}\footnote{The $X, Y$ here differ from $x, y$ in \cite{El} (2.2) by a minus-sign.}
$$
Y = \frac{Z_2}{ Z_3} \ , \ \ X = \frac{- Z_2^3}{ Z_3^2 Z_1 } \ .
$$
The automorphism $h$ of $K$ becomes $(X, Y) \mapsto (\frac{1}{1-X}, \frac{-X}{Y^3})$, and its  induced map on $\Gamma (K, \Omega)$ sends $\varphi_1, \varphi_2, \varphi_3$ to $\varphi_2, \varphi_3, \varphi_1$ respectively. Hence with the basis $\{ \varphi_1, \varphi_2, \varphi_3 \}$ of $\Gamma (K, \Omega)$, the $\langle g \rangle$-isomorphism $(i)$ in (\req(Z7)) is indeed a $\langle g, h \rangle$-isomorphism. 
It is known that one can embed  $K$ into ${\rm Jac}(K)$ via the integral of abelian differentials, 
\be
\iota \ ( = \iota_{p_0} ) \ : K \longrightarrow {\rm Jac}(K) \ , \ \ p \mapsto (\int_{p_0}^p \varphi_1 , \int_{p_0}^p \varphi_2 , \int_{p_0}^p \varphi_3 )^t \ , 
\ele(KimJ) 
where $p_0$ is an arbitrary fixed element in $K$. By choosing $p_0$ in $K^{ g }$, the embedding $\iota$ in (\req(KimJ)) is  $\langle g \rangle$-equivariant, by which a rational-curve  problem in CY 3-fold $\widehat{{\rm Jac}(K)/\kappa_7}$ was discussed in \cite{R03}. Furthermore through $\iota$, one can recover the Klein quartic $K$ in (\req(Kqrt)) as the canonical curve in the projective tangent space of ${\rm Jac}(K)$ at the identity element of ${\rm Jac}(K)$ (see, e.g. \cite{Mum}). Note that the group $\langle g, h \rangle$ has no common fixed point in $K$ since $h$ permutes elements in $K^g$ by sending $q_0, q_1, q_\infty$ to $q_1, q_\infty, q_0$ respectively. Hence there is no base element $p_0$ in the definition of $\iota$ in (\req(KimJ)) so that $K$ can be embedded into ${\rm Jac}(K)$ compatible with the $\langle g, h \rangle$-action.

Consider the quadratic imaginary subfield $\QZ(\eta)  \ (= \QZ ( \sqrt{-7}))$ of $\QZ (\mu)$, whose ring of integers $\ZZ[\mu]$ is a free $\ZZ[\eta]$-module of rank 3 with the basis $1, \mu, \mu^2$:
$$
\ZZ[\mu] = \ZZ[\eta] + \mu \ZZ[\eta] + \mu^2 \ZZ[\eta] \ .
$$
By the change of variables of $\CZ^3$ using the basis $\{ u(1), u(\mu), u(\mu^2) \}$, the transformation $m_\mu$ of $A( \mu )$ becomes (\req(eta7)). Hence we obtain the explicit description of $(ii)$ in (\req(Z7)). Furthermore, the automorphism $h$ of $A( \mu )$ corresponds to the following special automorphism $\widetilde{h}$ of $E_\eta^3$, 
$$
\widetilde{h} = \left( \begin{array}{ccc} 1 & \eta & 0 \\
0 & -1 &  1 \\
0 & -1 & 0 \end{array}
\right).
$$ 
Hence the isomorphism $(ii)$ in (\req(Z7)) can be further extended to the $\langle g, h \rangle$-equivariant one. 
By which and (\req(eta)), the $\widetilde{h}$-fixed point set of $E_\eta^3$ is the connected 1-torus $E_\eta \times [0]\times [0]$. 
Since one can write the generator $\alpha$ of $A( \mu )^{ m_\mu }$ in the form
$$\alpha = \frac{1}{7}\bigg(-u(1) - 5 u(\mu) -8 u(\mu^2) + 5 u(\eta) + 4 u(\mu \eta) -2 u(\mu^2 \eta) \bigg) ,
$$ 
the corresponding $g_\eta$-fixed point set in $E_\eta^3$ is the order 7 subgroup generated by $\frac{1}{7}([6+5 \eta], [2+ 4 \eta], [6+5 \eta])^t \in E_\eta^3$.

By choosing $p_0 \in K^{h}$ in (\req(KimJ)) to embed $K$ into ${\rm Jac}(K)$, through  isomorphisms in (\req(Z7)) and the projection of $E_\eta^3$ to its first factor $E_\eta$, $K$ becomes the cyclic $\langle h \rangle$-cover of $E_\eta$:
$$
K \ \longrightarrow \ K/\langle h \rangle \ \simeq \ E_\eta  \ .
$$
It is known that $E_\eta $ has the complex multiplication $\ZZ[\eta]$ with the following Weierstrass form:
$$
E_\eta : \ \ Y^2 = 4 X^3 + 21 X^2 + 28 X \ , \ \ (X, Y) \in \CZ^2 \ ,
$$ 
and one can write down an explicit relation between coordinates of $K$ and $E_\eta$ (see  \cite{El} \S 2).

\section{Intermediate Jacobian of Rigid CY 3-folds $\widehat{\bf T/G}$ }
First we give a simple derivation of 1-torus structure by its first cohomology group, a procedure which will be useful for later discussions of this section. 

For $\tau \in \HZ_+$, we write $\tau = \alpha + {\rm i} \beta$ with $\beta > 0$.
Denote by $z = x + {\rm i} y \in \CZ$ the uniformizing coordinates of the 1-torus $E_{\tau} $. The translations of $\CZ$ by $1$ and $\tau$ give rise to two basis elements of $H_1(E_{\tau}, \ZZ)$, denoted by $e, f$ respectively. Let $e^*, f^* $ be its dual base of  $H^1(E_{\tau}, \ZZ)$. We have the relations: 
$$
\left( \begin{array}{c} e \\ f \end{array}\right) = \left( \begin{array}{cc} 1 & 0 \\ \alpha & \beta   \end{array} \right) \left( \begin{array}{c} \partial_x \\  \partial_y  \end{array} \right) \ , \ \ \ \ (e^*, f^*) \left( \begin{array}{cc} 1 & 0 \\ \alpha & \beta \end{array} \right) = (d x , d y ) \ .
$$   
By the above second relation, one has 
$$
2 e^* = (1+ \frac{\alpha }{\beta} {\rm i}) d z + (1- \frac{\alpha }{\beta}{\rm i}) d \overline{z} \ , \ \ 2 f^* = \frac{1}{\beta }(-{\rm i}dz + {\rm i} d \overline{z})  \ .
$$
Hence through the base $\frac{1}{2} d \overline{z}$ of $H^{0,1}(E_{\tau}, \CZ)$, we can recover $E_\tau$ via the following isomorphisms:
\be
H^{0,1} (E_{\tau}, \CZ)/ H^1(E_{\tau}, \ZZ)  \ \simeq \ \CZ /(  \ZZ (1- \frac{\alpha }{\beta}{\rm i}) +  \ZZ \frac{\rm i }{\beta }) \, \ \simeq \, \ E_{\tau} \ .
\ele(Etau)

Now we consider the 3-torus $E_{\tau}^3$, with the uniformizing coordinates $(z_1, z_2, z_3) \in \CZ^3$, $z_j = x_j + {\rm i} y_j$ for $1 \leq j \leq 3$.  The integral basis of  $H_1(E_{\tau}^3, \ZZ)$ is given by $ \{ e_j, f_j \}_{j=1}^3 $, where $e_j, f_j$ are one-cycles of the $j$-th factor of $E_{\tau}^3$ for the coordinate $z_j$. Accordingly, $\{ e_j^*, f_j^* \}_{j=1}^3$ is its dual base of $H^1(E_{\tau}^3, \ZZ)$ so that the following relations hold:
$$
\left( \begin{array}{c} e_j \\ f_j \end{array}\right) = \left( \begin{array}{cc} 1 & 0 \\ \alpha & \beta   \end{array} \right) \left( \begin{array}{c} \partial_{x_j} \\  \partial_{y_j}  \end{array} \right) \ , \ \ \ \ (e_j^*, f_j^*) \left( \begin{array}{cc} 1 & 0 \\ \alpha & \beta \end{array} \right) = (d x_j , d y_j ) \ \ \ {\rm for } \ \ j=1, 2, 3 .
$$ 
Denote by $\Omega = d z_1 \wedge d z_2 \wedge d z_3$ the holomorphic volume form of $E_{\tau}^3$. One can express $\Omega$ in terms of $e_j^*$s, $f_j^*$s as follows:
\bea(cll)
\Omega &= & \Re (\Omega) +{\rm i} \Im (\Omega) \ ; \\
\Re (\Omega) &= & d x_1 \wedge d x_2 \wedge d x_3 - d x_1 \wedge d y_2 \wedge d y_3 - d y_1 \wedge d x_2 \wedge d y_3 - d y_1 \wedge d y_2 \wedge d x_3  \\
&= & e_1^* \wedge e_2^* \wedge e_3^* + \alpha ( e_1^* \wedge e_2^* \wedge f_3^* + e_1^* \wedge f_2^* \wedge e_3^* + f_1^* \wedge e_2^* \wedge e_3^* ) + \\
&& (\alpha^2 - \beta^2 ) (e_1^* \wedge f_2^* \wedge f_3^* + f_1^* \wedge e_2^* \wedge f_3^* + f_1^* \wedge f_2^* \wedge e_3^*) +( \alpha^3 - 3 \alpha \beta^2 ) f_1^* \wedge f_2^* \wedge f_3^* ;  \\
\Im (\Omega) &= & d x_1 \wedge d x_2 \wedge d y_3 + d x_1 \wedge d y_2 \wedge d x_3 + d y_1 \wedge d x_2 \wedge d x_3 - d y_1 \wedge d y_2 \wedge d y_3 \\
&= &  (3 \alpha^2 \beta - \beta^3  ) f_1^* \wedge f_2^* \wedge f_3^* + \beta ( e_1^* \wedge e_2^* \wedge f_3^* + e_1^* \wedge f_2^* \wedge e_3^* + f_1^* \wedge e_2^* \wedge e_3^* ) + \\
&& 2 \alpha \beta   (e_1^* \wedge f_2^* \wedge f_3^* + f_1^* \wedge e_2^* \wedge f_3^* + f_1^* \wedge f_2^* \wedge e_3^*)  \ .
\elea(vol)
In particular when $\tau = \omega ( = \frac{-1 + \sqrt{3} {\rm i}}{2})$, the volume form $\Omega =  \Re (\Omega) +{\rm i} \Im (\Omega)$ in (\req(vol)) has the following expression:
\be
\Re (\Omega)= A - \frac{1}{2} B  , \ 
\Im (\Omega) = \frac{\sqrt{3} }{2}  B  \ , \ \ \Leftrightarrow  \ \ (A , B) \left( \begin{array}{cc} 1 & 0 \\ - \frac{1}{2} & \frac{\sqrt{3}}{2} \end{array} \right) = (\Re (\Omega) , \Im (\Omega) )
\ele(vomega)
where $A, B$ are elements in $H^3(E_{\omega}^3, \ZZ)$ defined by
$$
\begin{array}{l}
A = e_1^* \wedge e_2^* \wedge e_3^* + f_1^* \wedge f_2^* \wedge f_3^* -e_1^* \wedge f_2^* \wedge f_3^* - f_1^* \wedge e_2^* \wedge f_3^* - f_1^* \wedge f_2^* \wedge e_3^* , \\
B = e_1^* \wedge e_2^* \wedge f_3^* + e_1^* \wedge f_2^* \wedge e_3^* + f_1^* \wedge e_2^* \wedge e_3^* - e_1^* \wedge f_2^* \wedge f_3^* - f_1^* \wedge e_2^* \wedge f_3^* - f_1^* \wedge f_2^* \wedge e_3^* \ .
\end{array}
$$
When $\tau = \eta ( = \frac{-1 + \sqrt{7} {\rm i}}{2})$, the volume form $\Omega$ in (\req(vol)) is expressed by
\be
\Re (\Omega) =  A - \frac{1}{2}  B , \ 
\Im (\Omega) = \frac{\sqrt{7}}{2} B , \ \ \Leftrightarrow  \ \ (A , B) \left( \begin{array}{cc} 1 & 0 \\ - \frac{1}{2} & \frac{\sqrt{7}}{2} \end{array} \right) = (\Re (\Omega) , \Im (\Omega) )
\ele(vnu)
where $A, B$ are elements in $H^3(E_{\eta}^3, \ZZ)$ defined by
$$
\begin{array}{ll}
A = & e_1^* \wedge e_2^* \wedge e_3^* + 2 (f_1^* \wedge f_2^* \wedge f_3^* -e_1^* \wedge f_2^* \wedge f_3^* - f_1^* \wedge e_2^* \wedge f_3^* - f_1^* \wedge f_2^* \wedge e_3^* ) , \\
B = &- f_1^* \wedge f_2^* \wedge f_3^* +  e_1^* \wedge e_2^* \wedge f_3^* + e_1^* \wedge f_2^* \wedge e_3^* + f_1^* \wedge e_2^* \wedge e_3^*  - e_1^* \wedge f_2^* \wedge f_3^* \\
&- f_1^* \wedge e_2^* \wedge f_3^* - f_1^* \wedge f_2^* \wedge e_3^* \ .
\end{array}
$$
Using (\req(vomega)) and (\req(vnu)), one can determine the intermediate Jacobian of rigid CY 3-folds $\widehat{T/G}$ as follows.  
\begin{theorem}\label{thm:InJ} 
For $(T, G) = (E_\omega^3, \langle m_\omega^3 \rangle), (E_\eta^3, \langle g_\eta \rangle)$, the intermediate Jacobian $J^2(\widehat{T/G})$ of the CY 3-fold $\widehat{T/G}$ is the 1-torus $E_\omega$, $E_\eta$ respectively, hence it is a CM elliptic curve defined over $\QZ$ . 
\end{theorem}
{\it Proof}. By the cohomology structure of the CY 3-folds $\widehat{T/G}$ and arithmetic properties of $E_\omega$ and $E_\eta$, it suffices to show that $H^3 (T,\ZZ)^G \ \simeq \ \ZZ^2$, and $H^{0,3}(T,\CZ)/H^3(T,\ZZ)^G$ is isomorphic to $E_\omega$ or $E_\eta$.  By the eigenvalues of a generator of $G$, $H^3(T,\CZ)^G$ is equal to $H^{3,0}(T,\CZ) + H^{0,3}(T,\CZ)$= $\CZ \Omega + \CZ \overline{\Omega}$. Hence $H^3(T,\RZ)^G$ is the 2-dimensional $\RZ$-space with the base $\{ \Re (\Omega), \Im (\Omega) \}$. Since $H^3(T,\ZZ)$ is a torsion-free $\ZZ$-module, $H^3 (T,\ZZ)^G$ is isomorphic to $\ZZ^2$. By the relations (\req(vomega)), (\req(vnu)), the expressions of $A, B$ in terms of $e_j^*, f_j^*$ imply that $A, B$ are elements in $H^3(T,\ZZ)^G$. Suppose  $H^3(T,\ZZ)^G = \ZZ C_1 \oplus \ZZ C_2$. Then one can express $A, B$ in the form 
$$
\left( \begin{array}{c} A \\ B \end{array} \right) = M \left( \begin{array}{c} C_1 \\ C_2 \end{array} \right)
$$
where $M$ is a square matrix of size 2 with integral entries. By comparing the coefficients of $e_1^* \wedge e_2^* \wedge e_3^*, e_1^* \wedge e_2^* \wedge f_3^*$ in the expressions of $A, B$ and $C_1, C_2$, one can conclude that $I_2 = M M'$ for some  2-by-2 integral matrix $M'$. Hence $M \in {\rm GL}_2(\ZZ)$, which implies $A, B$ form an integral basis of $H^3(T,\ZZ)^G$. By the relations (\req(vomega)), (\req(vnu)), when replacing $e^*, f^*, dx, dy$ by $A, B, \Re(\Omega), \Im (\Omega)$ in the argument of (\req(Etau)), one obtains the 1-torus structure of  $H^{0,3}(T,\CZ)/H^3(T,\ZZ)^G$, given by $E_\omega$, $E_\eta$ corresponding to $G = \ZZ_3, \ZZ_7$ respectively. 
$\Box$ \par

\section{$\QZ$-structure of Rigid CY 3-folds $\widehat{\bf T/G}$ }
In this section, we investigate the $\QZ$-structure of the Kummer-surface-type CY 3-folds $\widehat{T/G}$. One needs to identify the $\QZ$-structure of the orbifold $T/G$, then examine the rational toric structure of the crepant resolution $\widehat{ T/G}$. In general,
for a finite diagonal group $D$ action on $\CZ^n$, the orbifold $\CZ^n/D$ carries a natural toric structure (associated to the algebraic torus $\CZ^{* n}/D$). By methods in toric geometry, $\CZ^n/D$ is endowed with a $\QZ$-structure, so is and any toric variety over $\CZ^n/D$ with the birational $\QZ$-morphism. In particular, when $D = \langle \omega \rangle$ acting on $\CZ^3$, the $\QZ$-structures of $\CZ^3/\langle \omega \rangle$ and its crepant resolution are explicitly described as follows.
\begin{lemma} \label{lem:l3}
Let $U$ be the orbifold $\CZ^3/\langle \omega \rangle$, and $\widehat{U}$ the crepant (toric) resolution of $U$. Then $U$ and $\widehat{U}$ carry the canonical toric $\QZ$-structure so that the birational morphism of $\widehat{U}$ over $U$ is defined over $\QZ$. Furthermore, one has $U = {\rm Spec} \bigg( \CZ [\{ W_n \}_{n \in I}]/ J \bigg) $, where $I := \{n=(n_1, n_2, n_3) \in \ZZ_{\geq 0}^3, \sum_{i=1}^3 n_i = 3 \} $, and $J$ is the ideal generated by $W_n^3 - W_{(3,0,0)}^{n_1}W_{(0,3,0)}^{n_2}W_{(0,0,3)}^{n_3}$ for all $n \in I$.
\end{lemma}
{\it Proof}. In this proof, $(z_1, z_2, z_3)$ denote the coordinate system of $\CZ^3$. For $m= (m_1, m_2, m_3) \in \ZZ^3$, we denote the monomial $z_1^{m_1}z_2^{m_2}z_3^{m_3}$ by $z^m$. The $\langle \omega \rangle$-invariant polynomial subalgebra of $\CZ[z_1, z_2, z_3]$ is generated by $z^n$ for $n \in I$. By assigning $W_n$ to $z^n$, the structure of $U$ follows immediately. By its toric structure, the crepant resolution $\widehat{U}$ is composed of three affine charts ${\cal U}_k$'s with the following coordinate systems:
$$
\begin{array}{ll}
{\cal U}_1 \simeq \CZ^3 \ni (u_1, u_2, u_3) ,  & u_1 = \frac{z_1}{z_3} \ (= \frac{W_{(2,1,0)}}{W_{(1,1,1)}}) \ , u_2 = \frac{z_2}{z_3} \ (= \frac{W_{(1,2,0)}}{W_{(1,1,1)}})  \ ,u_3= z_3^3 \ (= W_{(0,0,3)}) \ ; \\
{\cal U}_2 \simeq \CZ^3 \ni (v_1, v_2, v_3) ,  & v_1 = \frac{z_2}{z_1} \ (= \frac{W_{(0,2,1)}}{W_{(1,1,1)}}) \ , v_2 = \frac{z_3}{z_1} \ (= \frac{W_{(0,1,2)}}{W_{(1,1,1)}})  \ ,v_3= z_1^3 \ (=W_{(3,0,0)}) \ ; \\
{\cal U}_3 \simeq \CZ^3 \ni (w_1, w_2, w_3) ,  & w_1 = \frac{z_1}{z_2} \ (= \frac{W_{(2,0,1)}}{W_{(1,1,1)}}) \ , w_2 = z_2^3 \ (= W_{(0, 3, 0)}) \ , \ w_3 = \frac{z_3}{z_2} (= \frac{W_{(1,0,2)}}{W_{(1,1,1)}})  \ .
\end{array}
$$
Hence the relations between affine coordinates of ${\cal U}_j$'s define the $\QZ$-structure of $\widehat{U}$. By $W_{(1,1,1)} = u_1u_2u_3 = v_1v_2v_3= w_1w_2w_3$, the expressions of affine coordinates in terms of $W_n$'s  describe the $\QZ$-morphism of $\widehat{U}$ over $U$.
$\Box$ \par \vspace{.2in}

As the quotient of $E_\omega$ by $\langle m_\omega \rangle $ is $\PZ^1$, $E_\omega$ can be represented as the smooth compactification of the following affine curve as in \cite{R03},
\be
Y^3 = X (X-1) \ , \ \ (X, Y ) \in \CZ^2 \ ,
\ele(cub)
with $dX/Y^2$ as the base element of holomorphic differentials, and the automorphism $m_\omega$ of $E_\omega$ corresponding to $(X, Y) \mapsto (X , \omega Y)$. Then the branched loci of the $X$-projection of (\req(cub)),  $X=0, 1, \infty$, corresponds three elements of $E_\omega^{ m_\omega }$, denoted by $e_0, e_1, e_\infty$ respectively. The involution of the torus $E_\omega$ becomes $(X, Y) \mapsto (\frac{1}{X}, \frac{-Y}{X})$, with  $e_1$ as the identity element of the order 3 $m_\omega$-fixed (additive) subgroup of $E_\omega$. Near $e_0$ and $e_1$, $Y$ provides the smooth coordinate of the curve (\req(cub)), while near $e_\infty$ the curve is represented by the variables, $\xi = 1/X , \eta =1/Y$, with the relation, $\xi^2 = \eta^3(1-\xi)$. Hence $t = \frac{\xi}{\eta} ( = \frac{Y}{X})$ serves as the smooth coordinate of (\req(cub)) near $e_\infty$, and $m_{\omega}$-transformation becomes, $t \mapsto \omega t$. 

We denote the $j$-th component of $E_\omega^3$ by $(X_j, Y_j)$ with the relation (\req(cub)) for $j=1,2, 3$. Now we describe the birational model of  $E_\omega^3/\langle m_\omega^3 \rangle$ using $(X_j, Y_j)$'s and the index set $I$ in Lemma \ref{lem:l3}. By identifying $Y_j$ with $z_j$ in Lemma \ref{lem:l3}, $E_\omega^3/\langle m_\omega^3 \rangle$ is birational to the affine variety defined by
$$
\begin{array}{l}
Y_j^3 = X_j (X_j-1)  \  \ , \ \ \  ( j=1, 2, 3 ) , \\ Y_1^3 = W_{(3,0,0)} , \ Y_2^3 = W_{(0,3,0)}  \ , Y_3^3 = W_{(0,0,3)} \ , \ \
W_n^3 = W_{(3,0,0)}^{n_1}W_{(0,3,0)}^{n_2}W_{(0,0,3)}^{n_3} \ \ ( n \in I ) \ ,
\end{array}
$$
which is equivalent to the following one in $\CZ^{|I|}$ defined by
$$
W_n^3 = \prod_{j=1}^3 X_j^{n_j} (X_j-1)^{n_j} \ \ {\rm for} \ n \in I \setminus \{(3,0, 0), (0,3,0),(0, 0, 3)\} \ .
$$
Indeed, the above affine variety is isomorphic to $\bigg( E_\omega^3/\langle m_\omega^3 \rangle \bigg) \setminus \bigcup_{j=1}^3 D_j $, where $D_j$ is the divisor in $E_\omega^3/\langle m_\omega^3 \rangle$ having $e_\infty$ as the $j$-th component in $E_\omega^3$. By  previous arithmetic discussions of the curve (\req(cub)) near $e_\infty$ and the crepant resolution in  Lemma \ref{lem:l3}, we obtain the following result. 
\begin{theorem}\label{thm:3eqn} 
For $(T, G) = (E_\omega^3, \langle m_\omega^3 \rangle )$, both $T/G$ and $\widehat{T/G}$ are defined over $\QZ$. Furthermore, $T/G$ is birational equivalent to ${\rm Spec} \bigg( \CZ [ X_1, X_2, X_3, \{ W_n \}_{n \in I^\prime }] / {\cal J} \bigg)$, where $I^\prime = I \setminus \{(3,0, 0), (0,3,0),(0, 0, 3)\}$ with $I$ in ${\rm Lemma} \ \ref{lem:l3}$, and ${\cal J}$ is the ideal generated by $W_n^3 - \prod_{j=1}^3 X_j^{n_j} (X_j-1)^{n_j}$ for $n \in I^\prime$.
\end{theorem}
$\Box$ \par \vspace{.2in}

Now we discuss the arithmetic structure of $\widehat{T/G}$ for $|G|=7$. 
It is known that the three-symmetric product ${\rm S}^3K$ of the curve $K$ is biregular to the blow-up of ${\rm Jac}(K)$ centered at an embedded curve $K$ \cite{Mum}. Indeed, by using the  embedding $\iota$ in (\req(KimJ)) with $p_0 \in K^{ g }$,  the birational morphism from ${\rm S}^3K$ onto ${\rm Jac}(K)$ is given by
\be
\phi: {\rm S}^3 K \ \longrightarrow \  {\rm Jac}(K) \ , \ \ \ p_1+ p_2+ p_3  \mapsto \sum_{j=1}^3 \iota (p_j) \ .
\ele(phi)
The automorphism $g^3 $ of $K^3$ induces an order 7 automorphism of ${\rm S}^3 K$, denoted by $\overline{g^3}$. Together with the action of $\kappa_7$ on ${\rm Jac}(K)$, the $\ZZ_7$-equivariant morphism $\phi$ induces a birational equivalence between ${\rm S}^3K/ \langle \overline{g^3} \rangle$ and ${\rm Jac}(K)/\kappa_7$.
We are going to use the affine curve (\req(Kaff)) to represent $K$, then derive the birational $\QZ$-structure of ${\rm Jac}(K)/\kappa_7$. First we note that $K^3$ is birational to the following affine variety in $\CZ^6$,
$$
Y_j^7 = X_j^2 (X_j-1) \ , \ \ j=1, 2, 3 \ , 
$$
whose $\QZ$-algebra structure is given by $A := \QZ[X_1, X_2, X_3, Y_1, Y_2, Y_3]/ J$ where $J$ is the ideal generated by $Y_j^7 - X_j^2 (X_j-1)$ for $j=1,2, 3$. Let $x_j, y_j$ be the generators of $A$ represented by $X_j, Y_j$ respectively, then $A= \QZ[x_1, x_2, x_3, y_1, y_2, y_3]$. The $\QZ$-algebra structure of ${\rm S}^3K$ is given by the subalgebra $S(A)$ of $A$.
$$
S(A) := \QZ [s_1, s_2, s_3, \sigma_1, \sigma_2, \sigma_3 ] \subset A \ 
$$
where $s_k, \sigma_k$ are the elementary symmetric functions of $x_j$'s, $y_j$'s respectively, i.e., $s_1 = \sum_{j=1}^3 x_j$, $s_2= \sum_{1 \leq i < j \leq 3} x_ix_j$, $s_3 = x_1x_2 x_3$, and similar expressions for $\sigma_k$ in terms of $y_j$'s. The action of $\overline{g^3}$ leaves $s_k$'s invariant, and sends $\sigma_k$ to $\mu^k \sigma_k$ for $ 1 \leq k \leq 3$. Hence the $\overline{g^3}$-invariant subalgebra $S(A)^{\langle \overline{g^3} \rangle}$ is the $\QZ$-subalgebra of $S(A)$ generated by
\be
s_1, \ s_2, \ s_3 , \ \sigma_1^7, \ \sigma_2^7, \ \sigma_3^7, \ \ \ \prod_{j=1}^3 \sigma_j^{n_j} \ \ \ {\rm for} \ 0 \leq n_j <7  \ \ {\rm with} \ n_1 + 2 n_2 + 3 n_3 \equiv 0 \pmod{7} \ . 
\ele(gT7q)
This implies that ${\rm S}^3K/ \langle \overline{g^3} \rangle$ is birational to ${\rm Spec} \bigg( S(A)^{\langle \overline{g^3} \rangle} \otimes_{\QZ} \CZ \bigg)$. Therefore we obtain the following result.
\begin{proposition} \label{prop:T7q}
For $(T, G) = ({\rm Jac}(K), \kappa_7 )$, $T/G$ is defined over $\QZ$. Indeed, $T/G$ (hence $\widehat{T/G}$) is birational equivalent to ${\rm Spec} \bigg( S(A)^{\langle \overline{g^3} \rangle} \otimes_{\QZ} \CZ \bigg)$, where $S(A)^{\langle \overline{g^3} \rangle}$ is the $\QZ$-subalgebra of $S(A)$ with the generators in (\req(gT7q)).
\end{proposition}
$\Box$ \par \noindent
{\bf Remark}. In the above proposition, we only established a birational $\QZ$-model of $T/G$ for $G \simeq \ZZ_7$. Even though the elements in (\req(gT7q)) form a set of generators for  the $\QZ$-algebra $S(A)^{\langle \overline{g^3} \rangle}$ representing $T/G$, a minimal set of generators and explicit relations among them remain interesting problems to be solved, and the answers will be useful for the quantitative study of arithmetic properties of $\widehat{T/G}$. 
$\Box$ \par \vspace{.2in}

We have introduced the $\QZ$-structure of ${\rm Jac}(K)$ from the rational structure of ${\rm S}^3K$ through the birational morphism $\phi$ in (\req(phi)). In order to further discuss the $\QZ$- structure of $\widehat{{\rm Jac}(K)/\kappa_7}$, one needs to examine the local structures induced  by ${\rm S}^3K$ under $g$-action near the $g$-fixed point set ${\rm Jac}(K)^g$. For this purpose, we need a more detailed relationship between ${\rm S}^3 K$ and ${\rm Jac}(K)$ under $\phi$.  Let $E$ be the exceptional divisor in ${\rm S}^3 K$ for $\phi$. Then $E$ is a $\PZ^1$-bundle over $\phi (E)$, which is a copy of $K$ in ${\rm Jac}(K)$. Indeed, the fiber-elements in $E$ over an element $x \in K \simeq \phi (E)$ can be described as follows. With $K$ in the form (\req(Kqrt)), for $x \in K$  we consider the linear system $E_x$ of (projective) lines in $\PZ^2$ passing through $x$. Then $E_x$ is a copy of $\PZ^1$, and for each $\ell \in E_x$, $\ell \cdot K$ is a degree 4 divisor in $K$ with its support containing $x$. Identify $\ell$ with the element $\ell \cdot K - x $ in  ${\rm S}^3K$, then $E_x$ gives rise to the fiber in $E$ over the element $x$. Denote by $\ell_{q_j}$ the line in $\PZ^2$ tangential to $K$ at $q_j \in K^g$ for $j=0, 1, \infty$. Then $ 
\ell_{q_0} \cdot K = 3 q_0 + q_1 ,  \ell_{q_1} \cdot K = 3 q_1 + q_\infty ,  \ell_{q_\infty} \cdot K = 3 q_\infty + q_0 $,
by which one has 
\be
\phi ( 2  q_0 + q_1 ) = \phi (3 q_\infty),  \ \ \phi ( 2  q_1 + q_\infty ) = \phi (3 q_0) , \ \  \phi ( 2  q_\infty + q_0 ) = \phi (3 q_1) .
\ele(phiq)

We now set $p_0= q_\infty$ in the definition of $\iota$ in (\req(KimJ)). Then $\phi$ defines an embedding of $K^2+ q_\infty$ into ${\rm Jac}(K)$, and it coincides with $\iota$ when restricted to $K+ 2q_\infty$.
Since $\langle g \rangle$ acts on $K$ freely outside $K^g$, the $\overline{g^3}$-fixed point set $({\rm S}^3K)^{\overline{g^3}}$ is equal to the three-symmetric product ${\rm S}^3(K^g)$ of $K^g$ with 10 elements. Under the morphism $\phi$, ${\rm S}^3(K^g)$ is mapped onto ${\rm Jac}(K)^g$ with (\req(phiq)) as the only multiple fibers. Other than those in (\req(phiq)),  the rest four elements of ${\rm Jac}(K)^g$ are in one-to-one correspondence with the following elements in ${\rm S}^3(K^g)$, 
\be 
q_1 + 2 q_\infty, \ \ \ 2q_0 + q_\infty, \ \ \ q_0+q_1+q_\infty, \ \ \ q_0+2q_1 \ \ \in {\rm S}^3(K^g) . 
\ele(JKg4)
Through the biregular isomorphism between ${\rm S}^3K \setminus E$ and ${\rm Jac}(K) \setminus \phi (E)$ via $\phi$, ${\rm Jac}(K) \setminus \phi (E)$  inherits the $\QZ$-structure of ${\rm S}^3K \setminus E$. We are going to examine the $g$-action near elements in (\req(JKg4)) to determine the $\QZ$-structure of $\widehat{{\rm Jac}(K)/\kappa_7}$ near  exceptional divisors associated to these $g$-fixed points in ${\rm Jac}(K)^g$. As before, $t_j$ denotes a local smooth coordinate system of the curve $K$ centered at $q_j \in K^g$ for $j=0, 1, \infty$. Furthermore, one may choose $t_j$ as one carrying a $\QZ$-structure form $K$. The local $\QZ$-structure of ${\rm S}^3K$ near an element in (\req(JKg4)) can be described by $t_j$'s in following manners : $(t_0, t_1, t_\infty)$ gives rise to a coordinate system near $q_0+q_1+q_\infty$, and $(t_i+t_i', t_it_i', t_j)$ is a smooth coordinates near $2q_i + q_j, (i \neq j)$ where $t_i'$ is another copy of $t_i$. By (\req(gloc)), these $\QZ$-coordinate systems are compatible with the matrix form of $g$ in (\req(ghr)), hence their corresponding crepant toric resolutions in $\widehat{{\rm Jac}(K)/\kappa_7}$ carry the canonical $\QZ$-structure. 

One can describe a $\QZ$-coordinate system of ${\rm Jac}(K)$ near an element in (\req(phiq)) from the $\QZ$-structure of  ${\rm S}^3K$ near the corresponding elements in $({\rm S}^3K)^{\overline{g^3}}$. For example, consider the element $0 = \phi ( 2  q_0 + q_1 ) = \phi (3 q_\infty)$. A $\QZ$-coordinate system near $2q_0 + q_1$ in ${\rm S}^3K$ is given by $(s_1, s_2, s_3):=(t_0+t_0', t_0t_0', t_1)$. Using relations in (\req(gloc)), one obtains a local expression of $\overline{g^3}$ in terms of $s_j$'s.
Since the $g$-action of ${\rm Jac}(K)$ near $0$ is equivalent to the matrix form of $g$ in (\req(ghr)), one concludes that a local coordinate system near $0$  in ${\rm Jac}(K)$, which is compatible with the $g$-action, is given by $(z_1, z_2, z_3):= (s_1, s_2s_3, s_3)$ or $(s_1, s_2s_3, s_2)$. Note that in this situation, $\phi(E)$ has the local defining equation $z_2=z_3=0$ in ${\rm Jac}(K)$, and $(z_1, \frac{z_2}{z_3})$ becomes a local coordinate system of $E \subset {\rm S}^3K $.  One may also determine the coordinate $z_j$'s in ${\rm Jac}(K)$  near $0$ by using the $\QZ$-coordinate system $s_j$'s near $3 q_\infty$ in ${\rm S}^3K$ with $(s_1, s_2, s_3)= (t_\infty + t_\infty' + t_\infty'',  t_\infty t_\infty' + t_\infty' t_\infty''+ t_\infty''t_\infty,  t_\infty t_\infty't_\infty'')$, then $(z_1, z_2, z_3) = (s_2, s_1, s_1s_3)$. 
Since $z_j$'s form a $\QZ$-coordinate system  compatible with the matrix form of $g$ in (\req(ghr)), the crepant toric resolution in $\widehat{{\rm Jac}(K)/\kappa_7}$ carries the canonical $\QZ$-structure near the exceptional divisors over $0$. By the same argument, one obtains the $\QZ$-structure of $\widehat{{\rm Jac}(K)/\kappa_7}$ near exceptional divisors associated to the other two elements in  (\req(phiq)). Therefore $\widehat{{\rm Jac}(K)/\kappa_7}$ has a $\QZ$-structure near all the exceptional divisors, which is compatible with the rational structure of ${\rm Jac}(K)/\kappa_7 $ in Proposition \ref{prop:T7q}. Hence we obtain following result.
\begin{theorem}\label{thm:7eqn} 
For $(T, G) = ({\rm Jac}(K), \kappa_7 )$, the CY $3$-fold $\widehat{T/G }$ is defined over $\QZ$. 
\end{theorem}
$\Box$ \par

\section{Modularity of Rigid CY 3-folds $\widehat{\bf T/G}$ } 
In this section, we study the modularity property of the Kummer-surface-type CY 3-folds $\widehat{T/G}$, (for the notion of modular CY varieties, see \cite{Yui}). By Theorems \ref{thm:3eqn} and \ref{thm:7eqn}, there is a projective CY 3-fold $V$ defined over $\QZ$ with $\overline{V} = \widehat{T/G}$. For a generic prime $\ell$, we consider the $\ell$-adic \'{e}tale 3rd cohomology group $H_{et}^3(\overline{V}, \QZ_\ell )$ of $V$, and its $L$-series $L(H_{et}^3(\overline{V}, \QZ_\ell ), s)$ is called the $L$-series of $V$. The modularity of $V$ guarantees the existence of a modular form $f$ so that the $L$-series of $V$ coincides with the $L$-series associated to $f$, up to a finite number of Euler factors. By Theorem \ref{thm:InJ}, the intermediate Jacobian $J^2(\overline{V})$ of $\overline{V}$ is a CM elliptic curve $E$ for an imaginary quadratic field $\KZ$ with $(E, \KZ)= (E_\omega, \QZ(\sqrt{-3})), (E_\eta , \QZ(\sqrt{-7})$ corresponding to $G= \ZZ_3, \ZZ_7$ respectively. By a classical result of Shimura and Taniyama \cite{ST}, $E$ is modular with its $L$-series given by a Hecke $L$-series for a Grossencharacter $\chi$, $L (E, s) = L (\chi, s )$. Conjecturally the $L$-series of $\overline{V}$ is given by 
$L(H_{et}^3(\overline{V}, \QZ_\ell ), s) = L( \chi^3, s)$ (see, Conjecture 8.4 in \cite{Yui}), hence $\overline{V}$ is modular. We now verify the modularity of $V$ from the arithmetic structure of its corresponding 3-torus $T$. 
Since the rational structures of $V$ and $T/G$ are derived from the $\QZ$-structure of $T$, we have  
$$
H_{et}^3(\overline{V}, \QZ_\ell ) \simeq H_{et}^3(T/G, \QZ_\ell ) \subset  H_{et}^3(T, \QZ_\ell ) \ .
$$
Hence $H_{et}^3(\overline{V}, \QZ_\ell )$ can be identified with a 2-dimensional irreducible factor of $H_{et}^3(T, \QZ_\ell )$ as modules of the Galois group ${\rm Gal}(\overline{\QZ}/\QZ)$. Indeed with an identification $H_{et}^3(T, \QZ_\ell )\simeq H^3(T, \CZ)$, $H_{et}^3(\overline{V}, \QZ_\ell )$ is the subspace $H^{3,0}(T) + H^{0,3}(T)$ of $H^3(T, \CZ)$. By the definition of $T$ and discussions in Sec. 2, the arithmetic structure of $T$ is isomorphic to the triple-product $E^3$ of $E$. It is known that there is a cusp form $g$ of weight 2 for some congruence subgroup $\Gamma_0(N)$ such that the $L$-series of $E$ is determined by $g$,
$$
L(E, s) = L(g, s) \sim \prod_{p: {\rm prime}} \frac{1}{1- a_p p^{-s} + p^{1-2s}} \ ,
$$
where $a_p$ (for a good prime $p$) is the coefficient in the characteristic polynomial of  Frobenius endomorphism on $H_{et}^1(E, \QZ_\ell )$ in the form $X^2 - a_p X + p = 0$, whose roots $\alpha_p, \beta_p$ have the absolute value $p^{\frac{1}{2}}$ with the argument given by the Grossencharacter $\chi$. With the identification $T = E^3$, one can see that the characteristic polynomial of Frobenius endomorphism on the $H_{et}^3(\overline{V}, \QZ_\ell )$-factor of $H_{et}^3(T, \QZ_\ell )$ has the roots $\alpha_p^3, \beta_p^3$, which implies $L(H_{et}^3(\overline{V}, \QZ_\ell ), s) = L( \chi^3, s)$. Hence we obtain
\begin{theorem}\label{thm:mod} 
For $(T, G) = (E_\omega^3, \langle m_\omega^3), ({\rm Jac}(K), \kappa_7 )$, the CY $3$-fold $\widehat{T/G }$ is modular. 
\end{theorem}
$\Box$ \par

\section{Concluding Remarks}
In this note, we revisit in details the structure of two Kummer-surface-type CY 3-folds $\widehat{T/G}$ with $|G| = 3, 7$ in \cite{RY}. The rigid property of these two CY 3-folds, whose effects could be aesthetically unsatisfactory for the 3-fold classification in algebraic geometry and likely irrelevant to the interest of string physics, has provided examples with a rich arithmetic structure due to the number-theoretic aspect connected to the construction. In this article, we have identified the intermediate Jacobian 1-tori of $\widehat{T/G}$, and provided the $\QZ$-structure of these CY 3-folds. The effort stems from    
conjectured modularity property for the class of rigid CY 3-folds \cite{Yui}. By the property of $T$ as the triple-product of a CM elliptic curve $E$, the modularity of $\widehat{T/G}$ follows from the classical result in \cite{ST} on modular nature of $E$. The conclusion is in accord with Conjecture 8.4 in \cite{Yui} about modularity of rigid CY 3-folds of CM type. 


\section*{Acknowledgements}
This article could be considered as a response to some questions raised in the talk given by the author in "BIRS Workshop: Calabi-Yau Varieties and Mirror Symmetry" December 6-11, 2003, in Banff Centre, Alberta, Canada. He would like to thank some in audience for the interesting remarks made during the talk. The author is grateful to N. Yui for information and correspondences on modularity of CY varieties related to this article. This work has been supported by NSC 92-2115-M-001-023.

\end{document}